\documentclass[12pt]{article}
\usepackage[margin=1in]{geometry}
\usepackage[T1]{fontenc}
\usepackage{amsmath,amssymb,amsthm,mathtools}
\usepackage{times}
\usepackage[colorlinks,linkcolor=blue,citecolor=blue,urlcolor=blue]{hyperref}
\hypersetup{
  pdftitle={Amplification of small L1 perturbations in the 2D Navier-Stokes equations near a self-similar vortex collapse},
  pdfauthor={Trinh T. Nguyen}}
\numberwithin{equation}{section}
\allowdisplaybreaks
\newtheorem{theorem}{Theorem}
\newtheorem{corollary}{Corollary}
\newtheorem{lemma}{Lemma}[section]
\newtheorem{proposition}[lemma]{Proposition}

\newcommand{\R}{\mathbb{R}}
\newcommand{\dd}{\,\mathrm d}
\DeclareMathOperator{\erf}{erf}
\newcommand{\gmm}{\gamma}
\newcommand{\epsz}{\varepsilon_0}
\newcommand{\Gtot}{\Gamma_{\mathrm{tot}}}
\newcommand{\Gabs}{\Gamma_{\mathrm{abs}}}
\newcommand{\omcol}{\omega_{\mathrm{col}}}
\newcommand{\omhat}{\widehat\omega}
\newcommand{\omtil}{\widetilde\omega}

\title{$L^1$ amplification near a self-similar vortex collapse \\
in 2D Navier-Stokes
}
\author{Trinh T. Nguyen\thanks{School of Mathematical Sciences,
Shanghai Jiao Tong University, Shanghai, China.\newline
Email: \texttt{q0004349@sjtu.edu.cn}.}}
\date{\today}

\begin{document}
\maketitle

\begin{abstract}
We establish quantitative $L^1$ instability for the two-dimensional
Navier--Stokes equations near a self-similar three-vortex collapse as viscosity tends to zero.
For every $0<\gmm<\frac12$ and $t_0\in(0,1)$, we construct
solutions that are $O(\nu^{1-\gmm})$-close in $L^1(\R^2)$
at $t_0$ but exhibit order-one $L^1$ separation in a shrinking
neighborhood of the collapse point at times approaching collapse
as $\nu\to0$.
By time translation, we also obtain pairs of solutions with
initial $L^1$ distance $O(\nu^{1-\gmm})$ that develop
order-one separation at a fixed positive time.
In both constructions, the solutions have the same viscosity
and total circulation and are uniformly bounded in $L^1(\R^2)$.
The proof exploits the amplification of small time delays, the fast speed of the vortex centers near collapse, together with
sharp inviscid approximation estimates for the Navier--Stokes
equations established in~\cite{Ga11,NguyenCollapse}.
\end{abstract}

\section{Introduction}

For viscous flows with small viscosity, a natural question is how small
an initial perturbation must be for the perturbed solution to remain close
to a given reference flow in a prescribed norm. One seeks to determine how
this size depends on viscosity and on the observation time, and at which
scales perturbations can produce an order-one separation. These stability
threshold questions arise more generally in evolution equations with small
physical parameters and have attracted considerable interest in
mathematical physics.
See
\cite{BHIW25,BMV16,BVW18,BCJ26,BGMZ25,CLN23,CLWZ20,CZG23,DG26,DGa26,DR04,Ga11,Gre00,GGN16,GN18,GN19,GRV24,KN24,KN25,MZ22,NguyenCollapse,WZ21,WZ23,ZZ25}
and the references therein.

In this work, we consider the two-dimensional Navier--Stokes
vorticity equation
\begin{equation}\label{eq:NS}
 \partial_t\omega+v^\omega\cdot\nabla\omega=\nu\Delta\omega,
 \qquad v^\omega=K*\omega,
 \qquad K(x)=\frac{x^\perp}{2\pi|x|^2},
\end{equation}
on $\R^2$, where $x^\perp=(-x_2,x_1)$ and $\nu>0$.

\paragraph{Quantitative $L^1$ stability of vortices.}
On fixed time intervals before point-vortex collision, Gallay
\cite{Ga11} establishes an $L^1$ approximation of the viscous
evolution by a superposition of Lamb--Oseen vortices, with
error $O(\nu)$. For a self-similar three-vortex collapse at
$t=1$ and every fixed $0<\gmm<1/2$, Nguyen \cite{NguyenCollapse}
extends this approximation up to time $1-\nu^{1-2\gmm}$,
with error $O_\gmm(\nu^{2\gmm})$.
As $\gmm\to0$, however, the error bound deteriorates and
ceases to be small when the distance to collapse becomes
comparable to the viscosity. Thus, the validity of the
inviscid approximation in this critical regime remains open.

 Fix $\Gamma\in\R$,
$M>|\Gamma|$, and $\nu_0>0$, and define
\[
 \mathcal X_{\Gamma,M}
 =\left\{f\in L^1(\R^2):
   \int_{\R^2}f(x)\dd x=\Gamma,\quad \|f\|_{L^1(\R^2)}\le M\right\}.
\]
We consider families $\Omega^\nu$ of global smooth solutions
of \eqref{eq:NS} with viscosity $\nu$ and initial vorticities
in $\mathcal X_{\Gamma,M}$, for every $0<\nu\le\nu_0$.
Conservation of circulation and the estimate
$\|\Omega^\nu(t)\|_{L^1(\R^2)}\le\|\Omega^\nu(0)\|_{L^1(\R^2)}$ imply that
$\Omega^\nu(t)\in\mathcal X_{\Gamma,M}$ for all $t\ge0$;
in particular,
\[
 \sup_{0<\nu\le\nu_0}\sup_{t\ge0}\|\Omega^\nu(t)\|_{L^1(\R^2)}\le M.
\]
It is therefore natural to ask for which positive scales
$r_\nu\to0$ there exist two
families $\Omega^\nu,\widetilde\Omega^\nu$ of global smooth
solutions of \eqref{eq:NS} with the same viscosity $\nu$,
initial vorticities
$\Omega^\nu(0),\widetilde\Omega^\nu(0)\in\mathcal X_{\Gamma,M}$,
and a time $T>0$, independent of $\nu$, such that
\begin{equation}\label{eq:threshold-question}
 \|\widetilde\Omega^\nu(0)-\Omega^\nu(0)\|_{L^1(\R^2)}\le r_\nu,
 \qquad
 \liminf_{\nu\to0}
 \|\widetilde\Omega^\nu(T)-\Omega^\nu(T)\|_{L^1(\R^2)}>0.
\end{equation}

Our main result, Theorem~\ref{thm:main}, establishes $L^1$ separation near
collapse. Let $\omcol^\nu$ denote the viscous
evolution of a fixed self-similar three-vortex collapsing
configuration. For every fixed $0<\gmm<1/2$
and $t_0\in(0,1)$, there is a solution $\omtil^\nu$
whose distance from $\omcol^\nu$ is at most
$\nu^{1-\gmm}$ at $t_0$ and at least a positive constant
independent of $\nu$ at $T_\nu=1-\nu^{1-2\gmm}$.
The separation occurs in a ball centered at the collapse point,
with radius of order $\nu^{1/2-\gmm}$.
Taking $\gmm>0$ arbitrarily small but fixed as $\nu\to0$
makes the exponent $1-\gmm$ in the initial perturbation size
arbitrarily close to $1$.

\paragraph{Notations.} 
The notation $A\lesssim B$ means that $A$ is bounded by
a constant multiple of $B$, independently of viscosity
and the variables in the estimate. Dependence on the fixed
collapsing configuration is suppressed; additional dependence
is indicated by subscripts. In particular, $\lesssim_I$ is
uniform for $t_0\in I$, and any dependence on $\nu$ is
displayed explicitly.
For positive quantities, $A_\nu\sim B_\nu$ means
$A_\nu/B_\nu\to1$ as $\nu\to0$.
The notation $o(1)$ denotes a quantity tending to zero as
$\nu\to0$, with all other parameters fixed. When $t_0$ ranges over a fixed compact
interval $I\subset(0,1)$, this convergence is uniform in $t_0$.
\section{The main result}\label{sec:main}
\paragraph{Self-similar vortex collapse solution \cite{Ar10,NguyenCollapse}.}
Fix nonzero real circulations $\Gamma_1,\Gamma_2,\Gamma_3$,
distinct points $z_1^0,z_2^0,z_3^0\in\R^2$, and
$\beta\in\R$ satisfying
\begin{equation}\label{eq:configuration}
 \frac12z_k^0-\beta(z_k^0)^\perp
 +\sum_{j\ne k}\Gamma_jK(z_k^0-z_j^0)=0,
 \qquad k=1,2,3.
\end{equation}
We call this a normalized collapsing configuration.
If $\mathfrak R_\theta$ denotes counterclockwise rotation
through angle $\theta$, the associated point-vortex paths are
\begin{equation}\label{eq:centers}
 z_k(t)=\sqrt{1-t}\,
 \mathfrak R_{-\beta\log(1-t)}z_k^0,
 \qquad 0\le t<1.
\end{equation}
They satisfy
\begin{equation}\label{eq:center-dynamics}
 \dot z_k(t)=\sum_{j\ne k}\Gamma_jK(z_k(t)-z_j(t)),
 \qquad
 |z_j(t)-z_k(t)|=\sqrt{1-t}\,|z_j^0-z_k^0|.
\end{equation}
Let $\omcol^\nu$ be the unique solution
of \eqref{eq:NS} with initial vorticity
\begin{equation}\label{eq:initial-measure}
 \omcol^\nu(0)
 =\sum_{k=1}^3\Gamma_k\delta_{z_k^0}.
\end{equation}
Write
\begin{equation}\label{eq:circulation-sums}
 \Gabs=\sum_{k=1}^3|\Gamma_k|,
 \qquad \Gtot=\sum_{k=1}^3\Gamma_k.
\end{equation}
Fix a normalized collapsing configuration, a compact interval
$I\subset(0,1)$, and a constant $C>\max_{1\le k\le3}|z_k^0|$.
For $0<\gmm<1/2$, set
\begin{equation}\label{eq:observation-radius}
 R_\nu=C\nu^{1/2-\gmm}.
\end{equation}
\begin{theorem}\label{thm:main}
There exists $c>0$, depending only on $I$ and the
configuration, such that for every fixed $0<\gmm<1/2$
and all sufficiently small $\nu>0$,
uniformly for $t_0\in I$, there is a solution $\omtil^\nu$
of \eqref{eq:NS} on $[t_0,\infty)$ with viscosity $\nu$ satisfying
\begin{equation}\label{eq:circulation}
 \int_{\R^2}\omtil^\nu(t,x)\dd x
 =\int_{\R^2}\omcol^\nu(t,x)\dd x
 =\Gtot,\qquad t\ge t_0,
\end{equation}
\begin{equation}\label{eq:mass-bound}
 \sup_{t\ge t_0}\bigl(
 \|\omtil^\nu(t)\|_{L^1(\R^2)}+
 \|\omcol^\nu(t)\|_{L^1(\R^2)}\bigr)
 \le2\Gabs,
\end{equation}
\begin{equation}\label{eq:initial-bound}
 \|\omtil^\nu(t_0)
       -\omcol^\nu(t_0)\|_{L^1(\R^2)}
 \le\nu^{1-\gmm},
\end{equation}
and, at $T_\nu:=1-\nu^{1-2\gmm}>\max I$,
\begin{equation}\label{eq:final-bound}
 \|\omtil^\nu(T_\nu)
       -\omcol^\nu(T_\nu)\|_{L^1(B(0,R_\nu))}\ge c.
\end{equation}
\end{theorem}
The viscosity threshold may also depend on $C$, while the constant
$c$ is independent of $C$ and $\gmm$.

The separation near collapse in Theorem~\ref{thm:main} yields the
following fixed-time consequence by time translation.

\begin{corollary}\label{cor:fixed-time}
With $R_\nu$ as in \eqref{eq:observation-radius}, there exists
$c>0$ such that, for every fixed $0<\gmm<1/2$ and all
sufficiently small $\nu>0$, there are solutions
$\Omega^\nu,\widetilde\Omega^\nu$ of \eqref{eq:NS} on
$[0,\infty)$ with viscosity $\nu$, satisfying
\[
 \int_{\R^2}\Omega^\nu(t,x)\dd x
 =\int_{\R^2}\widetilde\Omega^\nu(t,x)\dd x
 =\Gtot,\qquad t\ge0,
\]
\[
 \sup_{t\ge0}\bigl(
 \|\Omega^\nu(t)\|_{L^1(\R^2)}+
 \|\widetilde\Omega^\nu(t)\|_{L^1(\R^2)}\bigr)
 \le2\Gabs,
\]
and
\begin{equation}\label{eq:intro-result}
 \|\widetilde\Omega^\nu(0)-\Omega^\nu(0)\|_{L^1(\R^2)}
 \le\nu^{1-\gmm},
 \qquad
 \|\widetilde\Omega^\nu(\tfrac12)-\Omega^\nu(\tfrac12)\|_{L^1(B(0,R_\nu))}
 \ge c.
\end{equation}
\end{corollary}

\begin{proof}
First, we obtain solutions on $[0,\infty)$ with the required
circulation and $L^1$ bounds. Fix $I=[1/4,1/2]$ and set
\[
 s_\nu=T_\nu-\tfrac12
 =\tfrac12-\nu^{1-2\gmm}\in I
\]
for $\nu$ small enough that $\nu^{1-2\gmm}\le1/4$.
Apply Theorem~\ref{thm:main} with this interval $I$ and
$t_0=s_\nu$, for sufficiently small $\nu$, and define
\[
 \Omega^\nu(t)=\omcol^\nu(s_\nu+t),
 \qquad
 \widetilde\Omega^\nu(t)=\omtil^\nu(s_\nu+t),
 \qquad t\ge0.
\]
Time translation preserves \eqref{eq:NS}, and
\eqref{eq:circulation}--\eqref{eq:mass-bound} give the
circulation and $L^1$ bounds.

Next, we prove \eqref{eq:intro-result}. Since
$s_\nu+1/2=T_\nu$,
\[
 \begin{aligned}
 \|\widetilde\Omega^\nu(0)-\Omega^\nu(0)\|_{L^1(\R^2)}
 &=
 \|\omtil^\nu(s_\nu)
       -\omcol^\nu(s_\nu)\|_{L^1(\R^2)}
 \overset{\eqref{eq:initial-bound}}{\le}\nu^{1-\gmm},\\
 \|\widetilde\Omega^\nu(\tfrac12)-\Omega^\nu(\tfrac12)\|_{L^1(B(0,R_\nu))}
 &=
 \|\omtil^\nu(T_\nu)
       -\omcol^\nu(T_\nu)\|_{L^1(B(0,R_\nu))}
 \overset{\eqref{eq:final-bound}}{\ge}c.
 \end{aligned}
\]
The same constant $c$ and viscosity threshold apply for every
$s_\nu\in I$ by Theorem~\ref{thm:main}.
\end{proof}

\paragraph{Remark.}
Perturbations of size $O(\nu^{1-\gmm})$ in $L^1(\R^2)$ can
produce an order-one separation at $T_\nu=1-\nu^{1-2\gmm}$.
Describing the inviscid limit for this family all the way to
collapse time appears to be beyond current methods.
We refer to \cite{DT26,DM21,DMR24,DN19,GS26} and the
references therein for further discussions on this problem.

\section{Ideas of the proof}\label{sec:ideas}

\paragraph{$L^1$ separation of Gaussian profiles.}
We first describe the $L^1$ separation for a single vortex core.
For $\eta>0$, set
\[
 G(\xi)=\frac1{4\pi}e^{-|\xi|^2/4},
 \qquad
 P_z^\eta(x)=\frac1\eta
 G\!\left(\frac{x-z}{\sqrt\eta}\right).
\]
For two centers $z,\widetilde z\in\R^2$,
Lemma~\ref{lem:gaussian-translation} gives
\begin{equation}\label{eq:ideas-gaussian-distance}
 \|P_z^\nu-P_{\widetilde z}^\nu\|_{L^1(\R^2)}
 =2\erf\!\left(\frac{|\widetilde z-z|}{4\sqrt\nu}\right).
\end{equation}
Thus $|\widetilde z-z|\sim c\sqrt\nu$, with fixed $c>0$,
gives an order-one $L^1(\R^2)$ separation, whereas a displacement
$o(\sqrt\nu)$ gives a vanishing distance.
In our construction, the corresponding leading centers of the
two solution families agree at $t_0$. Near collapse, their gap
reaches a fixed positive fraction of the core width $\sqrt\nu$;
see \eqref{eq:ideas-core-gap} below. Formula
\eqref{eq:ideas-gaussian-distance} then detects this gap in $L^1$.

\paragraph{Tuning the circulation to slow down the vortices.}
Inspired by the time-rescaling observation in ~\cite{NC09},
we multiply all circulations by a common factor $0<\lambda_\nu<1$,
replacing $\Gamma_k$ by $\lambda_\nu\Gamma_k$.
By \eqref{eq:center-dynamics}, this multiplies the velocity of
every center at a fixed configuration by $\lambda_\nu$,
slowing down the vortex motion.

\paragraph{The corresponding change in vortex center motion.}
Fix $t_0\in(0,1)$. We compare two point-vortex systems
initialized at $t_0$ with identical center positions.
The original system has circulations $\Gamma_k$, while
the modified system has circulations $\lambda_\nu\Gamma_k$.
Reducing the circulations weakens the velocity field that
each vortex induces on the others. 
The centers follow the same trajectories more slowly,
so each moving vortex takes more time to reach a given
point on its trajectory.

This slower motion is described by
\begin{equation}\label{eq:ideas-clock}
 s_\nu(t)=t_0+\lambda_\nu(t-t_0),
 \qquad
 \widehat z_k(t)=z_k(s_\nu(t)).
\end{equation}
To verify that these paths correspond to the reduced
circulations, the chain rule gives
\[
 \begin{aligned}
 \dot{\widehat z}_k(t)
 &=\lambda_\nu\dot z_k(s_\nu(t))
 =\lambda_\nu\sum_{j\ne k}
   \Gamma_jK(z_k(s_\nu(t))-z_j(s_\nu(t)))\\
 &=\sum_{j\ne k}
   (\lambda_\nu\Gamma_j)K(\widehat z_k(t)-\widehat z_j(t)).
 \end{aligned}
\]
Hence $\widehat z_k$ are the point-vortex paths with
circulations $\lambda_\nu\Gamma_k$ and the same initial
positions $\widehat z_k(t_0)=z_k(t_0)$. At a later physical
time $t$, the modified system has reached the configuration
of the original system at the earlier time $s_\nu(t)<t$.
This is the meaning of slowing down the vortex motion.

\paragraph{Amplification near the vortex collapse.}
Fix $0<\gmm<1/2$ independently of $\nu$.
Choose $0<\epsz<1$ once and for all, depending only on the
fixed collapsing configuration, so that
\begin{equation}\label{eq:choice-epsilon}
 \epsz\bigl(3\Gabs+|\Gtot|\bigr)\le\frac12.
\end{equation}
Set
\begin{equation}\label{eq:ideas-parameters}
 \lambda_\nu=1-\epsz\nu^{1-\gmm},
 \qquad T_\nu=1-\nu^{1-2\gmm}.
\end{equation}
At the same physical time $T_\nu$, the original center is
at $z_k(T_\nu)$, while the modified center is at the earlier
point $z_k(s_\nu(T_\nu))$ on its trajectory.

For the self-similar paths \eqref{eq:centers},
\begin{equation}\label{eq:ideas-speed}
 |\dot z_k(t)|=\frac{|\dot z_k(0)|}{\sqrt{1-t}}.
\end{equation}
Each nonstationary center therefore moves faster near collapse.
For such a center, using $1-T_\nu=\nu^{1-2\gmm}$ and
$T_\nu\to1$, we formally obtain
\begin{equation}\label{eq:ideas-core-gap}
 \begin{aligned}
 |\widehat z_k(T_\nu)-z_k(T_\nu)|
 &\overset{\eqref{eq:ideas-clock}}{=}
 |z_k(s_\nu(T_\nu))-z_k(T_\nu)|\\
 &\overset{\eqref{eq:ideas-clock},\,\eqref{eq:center-shift}}{\sim}
 (1-\lambda_\nu)(T_\nu-t_0)|\dot z_k(T_\nu)|\\
 &\overset{\eqref{eq:ideas-parameters},\,\eqref{eq:ideas-speed}}{\sim}
 \epsz(1-t_0)|\dot z_k(0)|
   \frac{\nu^{1-\gmm}}{\sqrt{1-T_\nu}}
 \overset{\eqref{eq:ideas-parameters}}{\sim}
 \epsz(1-t_0)|\dot z_k(0)|\sqrt\nu.
 \end{aligned}
\end{equation}
This is precisely the positive center displacement relative to
the Gaussian width
required by \eqref{eq:ideas-gaussian-distance}.

\paragraph{Controlling the contributions of the other vortices.}
Fix a vortex whose center moves and, without loss of generality,
number it $1$. At time $T_\nu$, the two centers
$z_1(T_\nu)$ and $\widehat z_1(T_\nu)$ are $O(\sqrt\nu)$ apart
by \eqref{eq:ideas-core-gap}, so both lie in
$B(z_1(T_\nu),\nu^{(1-\gmm)/2})$ for small $\nu$.
The radius is much larger than the diffusion scale $\sqrt\nu$,
but much smaller than the distances to the other centers,
which have size $\nu^{1/2-\gmm}$ by \eqref{eq:center-dynamics}.
Proposition~\ref{prop:separation} shows that the other Gaussian
profiles contribute only $o(1)$ absolute vorticity inside this ball,
so they cannot cancel the positive $L^1$ separation between the
profiles centered at $z_1(T_\nu)$ and $\widehat z_1(T_\nu)$.

\paragraph{$L^1$ vorticity separation of Navier-Stokes.}
To transfer the Gaussian separation above to exact
Navier--Stokes solutions, we compare $\omcol^\nu$
with
\[
 \omhat^\nu(t,x)
 =\lambda_\nu\omcol^{\nu/\lambda_\nu}
   (s_\nu(t),x).
\]
This amplitude--time scaling reduces each leading circulation
to $\lambda_\nu\Gamma_k$ and moves its center along
$z_k(s_\nu(t))$, exactly as in the point-vortex calculation.
The transformed viscosity is
$\lambda_\nu(\nu/\lambda_\nu)=\nu$;
see Proposition~\ref{prop:comparison}.
At the fixed time $t_0$, the corresponding centers coincide.
Gallay's bound \eqref{eq:early} in
Proposition~\ref{prop:inviscid-approximation} and
Lemma~\ref{lem:gaussian-rescaling} therefore give
\[
 \|\omhat^\nu(t_0)
       -\omcol^\nu(t_0)\|_{L^1(\R^2)}
 \lesssim_{t_0}\nu+\epsz\nu^{1-\gmm}
 \sim\epsz\nu^{1-\gmm}.
\]
At $T_\nu=1-\nu^{1-2\gmm}$, Nguyen's bound
\eqref{eq:late}~\cite{NguyenCollapse}
gives a total approximation error $O_\gmm(\nu^{2\gmm})=o(1)$;
see \eqref{eq:two-approximations}.
The positive Gaussian $L^1$ separation therefore persists
for the exact Navier--Stokes solutions, as proved in
Proposition~\ref{prop:separation}.

\paragraph{Enforcing the total circulation constraint.}
The intermediate solution $\omhat^\nu$ has total
circulation $\lambda_\nu\Gtot$.
We place the missing circulation in a weak Gaussian centered
at $y_\nu\in\R^2$, far from the collapsing vortices, by prescribing
\[
 \omtil^\nu(t_0,x)
 =\omhat^\nu(t_0,x)
 +(1-\lambda_\nu)\Gtot G(x-y_\nu).
\]
Since $\int_{\R^2}G=1$, these data have the required total
circulation. The correction has $L^1(\R^2)$ norm
$(1-\lambda_\nu)|\Gtot|
=\epsz|\Gtot|\nu^{1-\gmm}$,
independently of $y_\nu$.
Moving the Gaussian farther away preserves its integral
and $L^1$ norm while weakening the velocity it induces
near the collapsing vortices, since $|K(x)|=1/(2\pi|x|)$.
For each fixed $\nu$, the interaction estimate \eqref{eq:decoupling} allows us to
choose $y_\nu$ sufficiently far away that the exact evolution
of these data differs from $\omhat^\nu$ by $o(1)$
in $L^1(B(0,R_\nu))$ at $T_\nu$, with $R_\nu$ as in
\eqref{eq:observation-radius}; see \eqref{eq:local-correction}.
The positive $L^1$ separation near collapse therefore persists.

\section{Proof of Theorem~\ref{thm:main}}\label{sec:proof}

Define
\begin{equation}\label{eq:gaussians}
 \begin{aligned}
 G(\xi)&=\frac1{4\pi}e^{-|\xi|^2/4},
 \qquad H_t^\nu(x)=\frac1{\nu t}
 G\!\left(\frac{x}{\sqrt{\nu t}}\right),\\
 Q^\nu(t,x)&=\sum_{k=1}^3\Gamma_kH_t^\nu(x-z_k(t)).
 \end{aligned}
\end{equation}
We use the Gaussian approximation at two times: at $t_0$
to estimate the initial perturbation, and near collapse
to detect the displacement of the vortex cores.

\begin{proposition}[\cite{Ga11,NguyenCollapse}]
\label{prop:inviscid-approximation}
Fix a normalized self-similar three-vortex collapsing
configuration as in \eqref{eq:configuration}.
Let $\omcol^\nu$ solve \eqref{eq:NS} with
initial data \eqref{eq:initial-measure}, and define $Q^\nu$
by \eqref{eq:gaussians}.
For every compact interval $I\subset(0,1)$ and all sufficiently
small $\nu>0$,
\begin{equation}\label{eq:early}
 \sup_{t\in I}
 \|\omcol^\nu(t)-Q^\nu(t)\|_{L^1(\R^2)}
 \lesssim_I\nu.
\end{equation}
For every fixed $0<\gmm<1/2$ and all sufficiently
small $\nu>0$,
\begin{equation}\label{eq:late}
 \sup_{0<t\le1-\nu^{1-2\gmm}}
 \|\omcol^\nu(t)-Q^\nu(t)\|_{L^1(\R^2)}
 \lesssim_\gmm\nu^{2\gmm}.
\end{equation}
\end{proposition}
Estimate \eqref{eq:early} follows from Theorem~2.4 and Lemma~2.3
of~\cite{Ga11}, including the replacement of the regularized centers
by $z_k(t)$. Estimate \eqref{eq:late} is Theorem~1.3
of~\cite{NguyenCollapse} with time exponent $1-2\gmm$.

\begin{lemma}\label{lem:gaussian-rescaling}
For every $\nu,t>0$ and $1/2\le\lambda_\nu<1$,
\begin{equation}\label{eq:gaussian-amplitude}
 \|\lambda_\nu H_t^{\nu/\lambda_\nu}-H_t^\nu\|_{L^1(\R^2)}
 \le3(1-\lambda_\nu).
\end{equation}
\end{lemma}

\begin{proof}
By \eqref{eq:gaussians},
\[
 \lambda_\nu H_t^{\nu/\lambda_\nu}(x)
 =\frac{\lambda_\nu^2}{\nu t}
 G\!\left(\frac{\sqrt{\lambda_\nu}\,x}{\sqrt{\nu t}}\right).
\]
The substitution $x=\sqrt{\nu t}\,\xi$, with
$\dd x=\nu t\,\dd\xi$, gives
\[
 \begin{aligned}
 &\|\lambda_\nu H_t^{\nu/\lambda_\nu}-H_t^\nu\|_{L^1(\R^2)}=
 \int_{\R^2}
 \left|\lambda_\nu^2G(\sqrt{\lambda_\nu}\,\xi)-G(\xi)\right|\dd\xi.
 \end{aligned}
\]
For $r\in[\lambda_\nu,1]$, define $F_r(\xi)=r^2G(\sqrt r\,\xi)$.
Since $\nabla G(\zeta)=-\zeta G(\zeta)/2$,
\[
 \begin{aligned}
 \partial_rF_r(\xi)
 &=2rG(\sqrt r\,\xi)
   +\frac{r^{3/2}}2\,\xi\cdot\nabla G(\sqrt r\,\xi)=r\left(2-\frac{r|\xi|^2}{4}\right)G(\sqrt r\,\xi).
 \end{aligned}
\]
This implies
\[
 \begin{aligned}
 \int_{\R^2}|\partial_rF_r(\xi)|\dd\xi
 &=\int_{\R^2}
   \left|2-\frac{|\zeta|^2}{4}\right|G(\zeta)\dd\zeta=\frac12\int_0^\infty
   \left|2-\frac{\rho^2}{4}\right|e^{-\rho^2/4}\rho\,\dd\rho\\
 &=\int_0^\infty|2-s|e^{-s}\dd s
 \le\int_0^\infty(2+s)e^{-s}\dd s=3.
 \end{aligned}
\]
This implies
\begin{equation}\label{eq:gaussian-rescaling}
 \begin{aligned}
 \|\lambda_\nu H_t^{\nu/\lambda_\nu}-H_t^\nu\|_{L^1(\R^2)}
 &=\int_{\R^2}
   \left|\int_{\lambda_\nu}^1\partial_rF_r(\xi)\dd r\right|\dd\xi\\
 &\le\int_{\lambda_\nu}^1\int_{\R^2}
   |\partial_rF_r(\xi)|\dd\xi\dd r
 \le3(1-\lambda_\nu).
 \end{aligned}
\end{equation}
\end{proof}

\begin{lemma}\label{lem:gaussian-translation}
For every $d\in\R^2$,
\begin{equation}\label{eq:gaussian-translation}
 \|G-G(\cdot-d)\|_{L^1(\R^2)}=2\erf(|d|/4),
 \qquad
 \erf(r)=\frac2{\sqrt\pi}\int_0^r e^{-s^2}\dd s.
\end{equation}
\end{lemma}

\begin{proof}
The identity holds for $d=0$. Suppose $d\ne0$.
By rotational invariance of $G$ and Lebesgue measure,
we may assume that $d=(|d|,0)$.
Then $G(\xi)>G(\xi-d)$ exactly when $\xi_1<|d|/2$.
Since both profiles have integral one, their difference has
integral zero, so its $L^1$ norm is twice the integral of
its positive part. Thus, we have
\[
 \begin{aligned}
 \|G-G(\cdot-d)\|_{L^1(\R^2)}
 &=2\int_{-\infty}^{|d|/2}\int_{\R}
   \bigl[G(s,u)-G(s-|d|,u)\bigr]\dd u\dd s\\
 &=\frac1{\sqrt\pi}\int_{-|d|/2}^{|d|/2}e^{-s^2/4}\dd s
 =\frac4{\sqrt\pi}\int_0^{|d|/4}e^{-r^2}\dd r
 =2\erf(|d|/4).
 \end{aligned}
\]
\end{proof}
\begin{proposition}\label{prop:comparison}
Fix a compact interval $I\subset(0,1)$ and take $\lambda_\nu$
as in \eqref{eq:ideas-parameters}.
For sufficiently small $\nu>0$ and $t_0\in I$, define
\begin{equation}\label{eq:intermediate}
 \omhat^\nu(t,x)
 =\lambda_\nu\omcol^{\nu/\lambda_\nu}
   \bigl(t_0+\lambda_\nu(t-t_0),x\bigr),
 \qquad t\ge t_0.
\end{equation}
This is a solution of \eqref{eq:NS} with viscosity $\nu$, and
\begin{equation}\label{eq:intermediate-mass}
 \|\omhat^\nu(t)\|_{L^1(\R^2)}
 \le\lambda_\nu\Gabs,\qquad
 \int_{\R^2}\omhat^\nu(t,x)\dd x
 =\lambda_\nu\Gtot.
\end{equation}
There is a constant $C_I$, independent of $\nu$, $\epsz$,
and $\gmm$, such that, uniformly for $t_0\in I$,
\begin{equation}\label{eq:comparison-initial}
 \|\omhat^\nu(t_0)
       -\omcol^\nu(t_0)\|_{L^1(\R^2)}
 \le C_I\nu^{1-\gmm}.
\end{equation}
\end{proposition}

\begin{proof}
First, we show that $\omhat^\nu$ satisfies
\eqref{eq:NS}. The Biot--Savart law is linear in the
vorticity, and the convolution acts only in the spatial variable.
For $t\ge t_0$ and $x\in\R^2$, \eqref{eq:intermediate} gives
\[
 \begin{aligned}
 v^{\omhat^\nu}(t,x)
 &=(K*\omhat^\nu(t,\cdot))(x)
 =\int_{\R^2}K(x-y)\omhat^\nu(t,y)\dd y\\
 &=\lambda_\nu\int_{\R^2}K(x-y)
   \omcol^{\nu/\lambda_\nu}
   (t_0+\lambda_\nu(t-t_0),y)\dd y
 =\lambda_\nu v^{\omcol^{\nu/\lambda_\nu}}
   (t_0+\lambda_\nu(t-t_0),x).
 \end{aligned}
\]
The chain rule therefore yields
\[
 \begin{aligned}
 &\partial_t\omhat^\nu
 +v^{\omhat^\nu}\cdot\nabla\omhat^\nu
 -\nu\Delta\omhat^\nu\\
 &\quad\overset{\eqref{eq:intermediate}}{=}
 \lambda_\nu^2\left[
 \partial_s\omcol^{\nu/\lambda_\nu}
 +v^{\omcol^{\nu/\lambda_\nu}}
       \cdot\nabla\omcol^{\nu/\lambda_\nu}
 -\frac{\nu}{\lambda_\nu}\Delta\omcol^{\nu/\lambda_\nu}
 \right]_{s=t_0+\lambda_\nu(t-t_0)}
 \overset{\eqref{eq:NS}}{=}0.
 \end{aligned}
\]

Next, we verify \eqref{eq:intermediate-mass}.
The advection--diffusion components with initial data
$\Gamma_k\delta_{z_k^0}$ preserve their signs and integrals.
Thus
\begin{equation}\label{eq:background-mass}
 \|\omcol^\nu(t)\|_{L^1(\R^2)}
 \le\Gabs,\qquad
 \int_{\R^2}\omcol^\nu(t,x)\dd x
 =\Gtot,\qquad t>0.
\end{equation}
Using \eqref{eq:intermediate}, we obtain
\[
 \begin{aligned}
 \|\omhat^\nu(t)\|_{L^1(\R^2)}
 &=\lambda_\nu
 \|\omcol^{\nu/\lambda_\nu}
       (t_0+\lambda_\nu(t-t_0))\|_{L^1(\R^2)}
 \overset{\eqref{eq:background-mass}}{\le}
 \lambda_\nu\Gabs,\\
 \int_{\R^2}\omhat^\nu(t,x)\dd x
 &=\lambda_\nu\int_{\R^2}
 \omcol^{\nu/\lambda_\nu}
       (t_0+\lambda_\nu(t-t_0),x)\dd x
 \overset{\eqref{eq:background-mass}}{=}
 \lambda_\nu\Gtot.
 \end{aligned}
\]

Finally, we prove \eqref{eq:comparison-initial}.
By \eqref{eq:ideas-parameters}, $\lambda_\nu\ge1/2$ for
sufficiently small $\nu$, so $\nu\le\nu/\lambda_\nu\le2\nu$.
For $t_0\in I$,
\begin{align}
 \|\omhat^\nu(t_0)
       -\omcol^\nu(t_0)\|_{L^1(\R^2)}
 &\overset{\eqref{eq:intermediate}}{=}
 \|\lambda_\nu\omcol^{\nu/\lambda_\nu}(t_0)
       -\omcol^\nu(t_0)\|_{L^1(\R^2)}\notag\\
 &\le
 \lambda_\nu\|\omcol^{\nu/\lambda_\nu}(t_0)
             -Q^{\nu/\lambda_\nu}(t_0)\|_{L^1(\R^2)}\notag\\
 &\quad+
 \|\omcol^\nu(t_0)-Q^\nu(t_0)\|_{L^1(\R^2)}
 +\|\lambda_\nu Q^{\nu/\lambda_\nu}(t_0)-Q^\nu(t_0)\|_{L^1(\R^2)}\notag\\
 &\overset{\eqref{eq:gaussians}}{\le}
 \lambda_\nu\|\omcol^{\nu/\lambda_\nu}(t_0)
             -Q^{\nu/\lambda_\nu}(t_0)\|_{L^1(\R^2)}\notag\\
 &\quad+
 \|\omcol^\nu(t_0)-Q^\nu(t_0)\|_{L^1(\R^2)}
 +\sum_k|\Gamma_k|
   \|\lambda_\nu H_{t_0}^{\nu/\lambda_\nu}-H_{t_0}^\nu\|_{L^1(\R^2)}\notag\\
 &\overset{\eqref{eq:gaussian-amplitude},\,\eqref{eq:ideas-parameters}}{\le}
 \lambda_\nu\|\omcol^{\nu/\lambda_\nu}(t_0)
             -Q^{\nu/\lambda_\nu}(t_0)\|_{L^1(\R^2)}\notag\\
 &\quad+
 \|\omcol^\nu(t_0)-Q^\nu(t_0)\|_{L^1(\R^2)}
 +3\epsz\Gabs\nu^{1-\gmm}
 \label{eq:comparison-initial-explicit}\\
 &\overset{\eqref{eq:early}}{\le}
 C_I\nu+3\epsz\Gabs\nu^{1-\gmm}
 \lesssim_I\nu^{1-\gmm}.
 \notag
\end{align}
\end{proof}

\begin{proposition}\label{prop:separation}
With $\epsz$ chosen in \eqref{eq:choice-epsilon}, fix
$0<\gmm<1/2$ and a compact interval $I\subset(0,1)$,
independently of $\nu$.
Set
\begin{equation}\label{eq:parameters}
 \lambda_\nu=1-\epsz\nu^{1-\gmm},
 \qquad T_\nu=1-\nu^{1-2\gmm}.
\end{equation}
Let $R_\nu$ be as in \eqref{eq:observation-radius}, and let
$\omhat^\nu$ be given by \eqref{eq:intermediate} for $t_0\in I$
with this choice of $\lambda_\nu$.
There exists $c_0>0$, depending only on $\epsz$, $I$,
and the fixed configuration, such that, for all sufficiently small $\nu$,
uniformly for $t_0\in I$,
\begin{equation}\label{eq:local-separation}
 \|\omhat^\nu(T_\nu)
       -\omcol^\nu(T_\nu)\|_{L^1(B(0,R_\nu))}
 \ge c_0,
\end{equation}
\end{proposition}

\begin{proof}
First, we show that the distance in \eqref{eq:local-separation}
can be replaced by the corresponding Gaussian distance
with a vanishing error. Set $h_\nu=(1-\lambda_\nu)(T_\nu-t_0)$.
This is the time lag accumulated by $\omhat^\nu$
between $t_0$ and $T_\nu$. For small $\nu$,
$\lambda_\nu\ge1/2$ and $T_\nu>\max I$, and
\eqref{eq:intermediate} gives
\begin{equation}\label{eq:lag}
 \omhat^\nu(T_\nu)
 =\lambda_\nu\omcol^{\nu/\lambda_\nu}
       (T_\nu-h_\nu),\qquad
 \frac{h_\nu}{\nu^{1-\gmm}}
 =\epsz(T_\nu-t_0),\qquad
 \frac{h_\nu}{\nu^{1-2\gmm}}\le \epsz\nu^\gmm.
\end{equation}
The modified centers in \eqref{eq:ideas-clock} satisfy
\begin{equation}\label{eq:shifted-centers}
 \widehat z_k(T_\nu)
 =z_k\bigl(t_0+\lambda_\nu(T_\nu-t_0)\bigr)
 =z_k(T_\nu-h_\nu).
\end{equation}
The Gaussian sums corresponding to
$\omcol^\nu(T_\nu)$ and
$\omhat^\nu(T_\nu)$ are, respectively,
\[
 \begin{aligned}
 Q^\nu(T_\nu,x)
 &=\sum_{k=1}^3\Gamma_kH_{T_\nu}^\nu(x-z_k(T_\nu)),\\
 \lambda_\nu Q^{\nu/\lambda_\nu}(T_\nu-h_\nu,x)
 &=\sum_{k=1}^3\lambda_\nu\Gamma_k
 H_{T_\nu-h_\nu}^{\nu/\lambda_\nu}(x-\widehat z_k(T_\nu)).
 \end{aligned}
\]
The second sum uses the viscosity, time, and amplitude factor
in \eqref{eq:lag}. We subtract the two profiles belonging to
each vortex and define, for $k=1,2,3$,
\begin{equation}\label{eq:gaussian-pair}
 f_{k,\nu}(x)=\Gamma_k\left[
 H_{T_\nu}^\nu(x-z_k(T_\nu))
 -\lambda_\nu H_{T_\nu-h_\nu}^{\nu/\lambda_\nu}
       (x-\widehat z_k(T_\nu))
 \right].
\end{equation}
Thus the difference of the two Gaussian sums is exactly
\[
 \sum_{k=1}^3f_{k,\nu}
 =Q^\nu(T_\nu)-\lambda_\nu Q^{\nu/\lambda_\nu}(T_\nu-h_\nu).
\]
Subtracting this identity from the difference of the exact
solutions leaves only their Gaussian approximation errors:
\[
 \begin{aligned}
 &\omcol^\nu(T_\nu)-\omhat^\nu(T_\nu)
       -\sum_{k=1}^3f_{k,\nu}\\
 &\quad=\omcol^\nu(T_\nu)-Q^\nu(T_\nu)
 -\lambda_\nu\left[
 \omcol^{\nu/\lambda_\nu}(T_\nu-h_\nu)
       -Q^{\nu/\lambda_\nu}(T_\nu-h_\nu)\right].
 \end{aligned}
\]
Reversing the sign does not change an $L^1$ norm. The reverse
triangle inequality, followed by the triangle inequality and
$B(0,R_\nu)\subset\R^2$, therefore gives
\begin{equation}\label{eq:gaussian-reduction}
 \begin{aligned}
 &\left|
 \|\omhat^\nu(T_\nu)-\omcol^\nu(T_\nu)
       \|_{L^1(B(0,R_\nu))}
 -\left\|\sum_k f_{k,\nu}\right\|_{L^1(B(0,R_\nu))}
 \right|\\
 &\quad\le
 \left\|\omcol^\nu(T_\nu)-\omhat^\nu(T_\nu)
       -\sum_{k=1}^3f_{k,\nu}\right\|_{L^1(B(0,R_\nu))}\\
 &\quad\le
 \lambda_\nu
 \|\omcol^{\nu/\lambda_\nu}(T_\nu-h_\nu)
       -Q^{\nu/\lambda_\nu}(T_\nu-h_\nu)\|_{L^1(\R^2)}
 +\|\omcol^\nu(T_\nu)-Q^\nu(T_\nu)\|_{L^1(\R^2)}.
 \end{aligned}
\end{equation}
Now we check the time condition needed to apply 
\eqref{eq:late} to the first error in
\eqref{eq:gaussian-reduction}. At viscosity $\nu/\lambda_\nu$,
we need
\begin{equation}\label{eq:shifted-approximation-time}
 0<T_\nu-h_\nu
 \le1-\left(\frac{\nu}{\lambda_\nu}\right)^{1-2\gmm}.
\end{equation}
The lower bound follows from
$T_\nu-h_\nu=t_0+\lambda_\nu(T_\nu-t_0)>0$.
Since $T_\nu=1-\nu^{1-2\gmm}$, the upper bound is equivalent to
\[
 (\nu/\lambda_\nu)^{1-2\gmm}-\nu^{1-2\gmm}\le h_\nu.
\]
We first factor out $\nu^{1-2\gmm}$ and then apply the
fundamental theorem of calculus:
\[
 \begin{aligned}
 (\nu/\lambda_\nu)^{1-2\gmm}-\nu^{1-2\gmm}
 &=\nu^{1-2\gmm}\bigl(\lambda_\nu^{-(1-2\gmm)}-1\bigr)=\nu^{1-2\gmm}(1-2\gmm)
   \int_{\lambda_\nu}^1 r^{2\gmm-2}\dd r.
 \end{aligned}
\]
Since $0<\gmm<1/2$ and $1/2\le\lambda_\nu\le r\le1$,
\[
 r^{2\gmm-2}=r^{2\gmm}r^{-2}\le r^{-2}\le4.
\]
Integrating this bound over $[\lambda_\nu,1]$ yields
\[
 \int_{\lambda_\nu}^1 r^{2\gmm-2}\dd r
 \le\int_{\lambda_\nu}^1 4\dd r=4(1-\lambda_\nu).
\]
Substituting this bound into the identity above and using
$0<1-2\gmm<1$, we obtain
\[
 0\le(\nu/\lambda_\nu)^{1-2\gmm}-\nu^{1-2\gmm}
 \le4\nu^{1-2\gmm}(1-2\gmm)(1-\lambda_\nu)
 \le4\nu^{1-2\gmm}(1-\lambda_\nu).
\]
We now divide by $h_\nu=(1-\lambda_\nu)(T_\nu-t_0)>0$.
The factor $1-\lambda_\nu$ cancels, giving
\[
 0\le
 \frac{(\nu/\lambda_\nu)^{1-2\gmm}-\nu^{1-2\gmm}}{h_\nu}
 \le\frac{4\nu^{1-2\gmm}(1-\lambda_\nu)}
          {(1-\lambda_\nu)(T_\nu-t_0)}
 =\frac{4\nu^{1-2\gmm}}{T_\nu-t_0}.
\]
Since $t_0<1$, we may take $\nu$ sufficiently small that
$\nu^{1-2\gmm}\le(1-t_0)/2$. Then
\[
 T_\nu-t_0=1-t_0-\nu^{1-2\gmm}
 \ge\frac{1-t_0}{2}>0.
\]
Substituting this lower bound into the preceding estimate gives
\[
 0\le
 \frac{(\nu/\lambda_\nu)^{1-2\gmm}-\nu^{1-2\gmm}}{h_\nu}
 \le\frac{8}{1-t_0}\nu^{1-2\gmm}\longrightarrow0,
\]
This proves \eqref{eq:shifted-approximation-time}, uniformly when
$t_0$ stays bounded away from $1$.
We may therefore apply \eqref{eq:late} at viscosity
$\nu/\lambda_\nu$ and time $T_\nu-h_\nu$, as well as at
viscosity $\nu$ and time $T_\nu$. Since $\nu/\lambda_\nu\le2\nu$,
we obtain
\begin{equation}\label{eq:two-approximations}
 \begin{aligned}
 &\lambda_\nu
 \|\omcol^{\nu/\lambda_\nu}(T_\nu-h_\nu)
       -Q^{\nu/\lambda_\nu}(T_\nu-h_\nu)\|_{L^1(\R^2)}\\
 &\quad+
 \|\omcol^\nu(T_\nu)-Q^\nu(T_\nu)\|_{L^1(\R^2)}
 \overset{\eqref{eq:late}}{\lesssim}_\gmm\nu^{2\gmm}.
 \end{aligned}
\end{equation}

Combining \eqref{eq:gaussian-reduction} and
\eqref{eq:two-approximations}, we obtain
\begin{equation}\label{eq:gaussian-lower-bound}
 \|\omhat^\nu(T_\nu)
       -\omcol^\nu(T_\nu)\|_{L^1(B(0,R_\nu))}
 \ge\left\|\sum_{k=1}^3 f_{k,\nu}\right\|_{L^1(B(0,R_\nu))}
       -C_\gmm\nu^{2\gmm},
\end{equation}
where $C_\gmm$ is independent of $\nu$ and $t_0$.
Since $\nu^{2\gmm}\to0$, it remains to bound the Gaussian
sum from below by a positive constant independent of $\nu$.
For each $k=1,2,3$, write
\begin{equation}\label{eq:gaussian-pair-decomposition}
 \begin{aligned}
 f_{k,\nu}(x)
 &=\underbrace{\Gamma_k\left[
 H_{T_\nu}^\nu(x-z_k(T_\nu))
 -H_{T_\nu}^\nu(x-\widehat z_k(T_\nu))
 \right]}_{\mathrm{(I)}}\\
 &\quad+\underbrace{\Gamma_k\left[
 H_{T_\nu}^\nu(x-\widehat z_k(T_\nu))
 -\lambda_\nu H_{T_\nu-h_\nu}^{\nu/\lambda_\nu}
       (x-\widehat z_k(T_\nu))
 \right]}_{\mathrm{(II)}}.
 \end{aligned}
\end{equation}
We first estimate $\mathrm{(II)}$.
Translation invariance gives
\[
 \|\mathrm{(II)}\|_{L^1(\R^2)}
 =|\Gamma_k|\,
 \|H_{T_\nu}^\nu-\lambda_\nu H_{T_\nu-h_\nu}^{\nu/\lambda_\nu}
       \|_{L^1(\R^2)}.
\]
We write
\[
 \begin{aligned}
 H_{T_\nu}^\nu-\lambda_\nu H_{T_\nu-h_\nu}^{\nu/\lambda_\nu}
 &=\left(H_{T_\nu-h_\nu}^\nu
       -\lambda_\nu H_{T_\nu-h_\nu}^{\nu/\lambda_\nu}\right)\\
 &\quad+\left(\frac{T_\nu-h_\nu}{T_\nu}H_{T_\nu}^\nu
       -H_{T_\nu-h_\nu}^\nu\right)
 +\frac{h_\nu}{T_\nu}H_{T_\nu}^\nu.
 \end{aligned}
\]
The triangle inequality and $\|H_{T_\nu}^\nu\|_{L^1(\R^2)}=1$ yield
\[
 \begin{aligned}
 \|\mathrm{(II)}\|_{L^1(\R^2)}
 &\le|\Gamma_k|\Biggl(
 \|H_{T_\nu-h_\nu}^\nu
       -\lambda_\nu H_{T_\nu-h_\nu}^{\nu/\lambda_\nu}\|_{L^1(\R^2)}+
 \left\|\frac{T_\nu-h_\nu}{T_\nu}H_{T_\nu}^\nu
       -H_{T_\nu-h_\nu}^\nu\right\|_{L^1(\R^2)}
 +\frac{h_\nu}{T_\nu}\Biggr).
 \end{aligned}
\]
Taking $t=T_\nu-h_\nu$ in Lemma~\ref{lem:gaussian-rescaling} gives
\[
 \|\lambda_\nu H_{T_\nu-h_\nu}^{\nu/\lambda_\nu}
       -H_{T_\nu-h_\nu}^\nu\|_{L^1(\R^2)}
 \le3(1-\lambda_\nu).
\]
Next, take $t=T_\nu-h_\nu$ and replace $\lambda_\nu$ in
Lemma~\ref{lem:gaussian-rescaling} by $(T_\nu-h_\nu)/T_\nu\in[1/2,1)$.
Since \eqref{eq:gaussians} gives
$H_{T_\nu-h_\nu}^{\nu T_\nu/(T_\nu-h_\nu)}=H_{T_\nu}^\nu$,
we obtain
\[
 \left\|\frac{T_\nu-h_\nu}{T_\nu}H_{T_\nu}^\nu
       -H_{T_\nu-h_\nu}^\nu\right\|_{L^1(\R^2)}
 \le3\left(1-\frac{T_\nu-h_\nu}{T_\nu}\right)
 =\frac{3h_\nu}{T_\nu}.
\]
Substituting these two bounds, we obtain
\[
 \|\mathrm{(II)}\|_{L^1(\R^2)}
 \le|\Gamma_k|\left(3(1-\lambda_\nu)+\frac{4h_\nu}{T_\nu}\right)
 \le7|\Gamma_k|\epsz\nu^{1-\gmm}\longrightarrow0,
\]
where, since $0<t_0<T_\nu$,
\[
 \frac{h_\nu}{T_\nu}
 \overset{\eqref{eq:lag}}{=}
 (1-\lambda_\nu)\frac{T_\nu-t_0}{T_\nu}
 \le1-\lambda_\nu
 \overset{\eqref{eq:parameters}}{=}\epsz\nu^{1-\gmm}.
\]

For $\mathrm{(I)}$, the change of variables
$x=z_k(T_\nu)+\sqrt{\nu T_\nu}\,\xi$ and
Lemma~\ref{lem:gaussian-translation} give
\[
 \|\mathrm{(I)}\|_{L^1(\R^2)}
 =2|\Gamma_k|\erf\!\left(
   \frac{|\widehat z_k(T_\nu)-z_k(T_\nu)|}{4\sqrt{\nu T_\nu}}
 \right).
\]
The reverse triangle inequality yields
\begin{equation}\label{eq:pair-width-error}
 \|f_{k,\nu}\|_{L^1(\R^2)}
 \ge\|\mathrm{(I)}\|_{L^1(\R^2)}-O(\nu^{1-\gmm}).
\end{equation}
By \eqref{eq:shifted-centers} and the fundamental theorem of calculus,
\[
 \begin{aligned}
 d_{k,\nu}:=\frac{\widehat z_k(T_\nu)-z_k(T_\nu)}{\sqrt\nu}
 &=\frac{z_k(T_\nu-h_\nu)-z_k(T_\nu)}{\sqrt\nu}\\
 &=-\frac1{\sqrt\nu}\int_{T_\nu-h_\nu}^{T_\nu}\dot z_k(t)\dd t.
 \end{aligned}
\]
Applying it once more gives
\[
 \begin{aligned}
 d_{k,\nu}+\frac{h_\nu}{\sqrt\nu}\dot z_k(T_\nu)
 &=\frac1{\sqrt\nu}\int_{T_\nu-h_\nu}^{T_\nu}
   \bigl(\dot z_k(T_\nu)-\dot z_k(t)\bigr)\dd t\\
 &=\frac1{\sqrt\nu}\int_{T_\nu-h_\nu}^{T_\nu}
   \int_t^{T_\nu}\ddot z_k(r)\dd r\dd t.
 \end{aligned}
\]
Set $b_k=\frac12z_k^0-\beta(z_k^0)^\perp$.
Differentiating \eqref{eq:centers} twice gives
\begin{equation}\label{eq:center-derivatives}
 \ddot z_k(t)
 =-\frac{\mathfrak R_{-\beta\log(1-t)}
       (\frac12 b_k+\beta b_k^\perp)}{(1-t)^{3/2}}.
\end{equation}
Since rotations preserve length,
\[
 |\ddot z_k(t)|
 =\frac{|\frac12 b_k+\beta b_k^\perp|}{(1-t)^{3/2}}
 \lesssim\frac1{(1-t)^{3/2}}.
\]
Substituting this bound into the double integral above and using
$1-r\ge1-T_\nu$ for $r\le T_\nu$ gives
\[
 \begin{aligned}
 \left|d_{k,\nu}+\frac{h_\nu}{\sqrt\nu}\dot z_k(T_\nu)\right|
 &\le\frac1{\sqrt\nu}\int_{T_\nu-h_\nu}^{T_\nu}
       \int_t^{T_\nu}|\ddot z_k(r)|\dd r\dd t\\
 &\lesssim\frac1{\sqrt\nu}\int_{T_\nu-h_\nu}^{T_\nu}
       \int_t^{T_\nu}\frac{\dd r\dd t}{(1-r)^{3/2}}\\
 &\le\frac1{\sqrt\nu(1-T_\nu)^{3/2}}
       \int_{T_\nu-h_\nu}^{T_\nu}(T_\nu-t)\dd t
 =\frac{h_\nu^2}{2\sqrt\nu(1-T_\nu)^{3/2}}\\
 &\overset{\eqref{eq:parameters},\,\eqref{eq:lag}}{=}
 \frac{\epsz^2\nu^{2-2\gmm}(T_\nu-t_0)^2}
      {2\nu^{1/2}(\nu^{1-2\gmm})^{3/2}}
 =\frac{\epsz^2}{2}(T_\nu-t_0)^2\nu^\gmm
 \le\frac{\epsz^2}{2}\nu^\gmm.
 \end{aligned}
\]
Here $0<T_\nu-t_0<1$.
Moreover, \eqref{eq:centers}, \eqref{eq:parameters},
and \eqref{eq:lag} give
\[
 \frac{h_\nu}{\sqrt\nu}|\dot z_k(T_\nu)|
 =\frac{h_\nu|b_k|}{\sqrt\nu\sqrt{1-T_\nu}}
 =\frac{h_\nu|b_k|}{\nu^{1-\gmm}}
 =\epsz(T_\nu-t_0)|b_k|.
\]
The reverse triangle inequality therefore yields
\begin{equation}\label{eq:center-shift}
 \left||d_{k,\nu}|-\epsz(T_\nu-t_0)|b_k|\right|
 \lesssim_{\epsz}\nu^\gmm,
 \qquad
 |d_{k,\nu}|\longrightarrow\epsz(1-t_0)|b_k|.
\end{equation}
Combining \eqref{eq:pair-width-error} with \eqref{eq:center-shift}
and $T_\nu\to1$, we obtain
\begin{equation}\label{eq:pair-limit}
 \|f_{k,\nu}\|_{L^1(\R^2)}
 \ge2|\Gamma_k|\erf\!\left(
   \frac{\epsz(1-t_0)|b_k|}{4}\right)-o(1).
\end{equation}

To conclude \eqref{eq:local-separation} from
\eqref{eq:gaussian-lower-bound} and \eqref{eq:pair-limit},
it suffices to prove
\[
 \left\|\sum_{k=1}^3 f_{k,\nu}\right\|_{L^1(B(0,R_\nu))}
 \ge\sum_{k=1}^3\|f_{k,\nu}\|_{L^1(\R^2)}-o(1).
\]
For each vortex $k$, let
$B_{k,\nu}=B(z_k(T_\nu),\nu^{(1-\gmm)/2})$.
We decompose
\[
 f_{k,\nu}
 =f_{k,\nu}\mathbf1_{B_{k,\nu}}
 +f_{k,\nu}\mathbf1_{\R^2\setminus B_{k,\nu}},
\]
where $\mathbf1_A$ denotes the indicator of $A$.
We now prove that
\[
 B_{i,\nu}\cap B_{j,\nu}=\varnothing\quad(i\ne j),
 \qquad B_{k,\nu}\subset B(0,R_\nu).
\]
For $i\ne j$, \eqref{eq:center-dynamics} and
$1-T_\nu=\nu^{1-2\gmm}$ give
\[
 \begin{aligned}
 |z_i(T_\nu)-z_j(T_\nu)|
 &=\sqrt{1-T_\nu}\,|z_i^0-z_j^0|\\
 &=(\nu^{1-2\gmm})^{1/2}|z_i^0-z_j^0|
 =\nu^{1/2-\gmm}|z_i^0-z_j^0|.
 \end{aligned}
\]
Choose $\nu$ sufficiently small that
$2\nu^{\gmm/2}<|z_i^0-z_j^0|$ for every $i\ne j$.
Then
\begin{equation}\label{eq:ball-center-distance}
 \begin{aligned}
 |z_i(T_\nu)-z_j(T_\nu)|
 &=\nu^{1/2-\gmm}|z_i^0-z_j^0|
 >2\nu^{1/2-\gmm}\nu^{\gmm/2}\\
 &=2\nu^{1/2-\gmm+\gmm/2}
 =2\nu^{(1-\gmm)/2}.
 \end{aligned}
\end{equation}
If there exists an $x\in B_{i,\nu}\cap B_{j,\nu}$ with $i\ne j$, then
\[
 \begin{aligned}
 |z_i(T_\nu)-z_j(T_\nu)|
 &\le |z_i(T_\nu)-x|+|x-z_j(T_\nu)|<2\nu^{(1-\gmm)/2},
 \end{aligned}
\]
contradicting \eqref{eq:ball-center-distance}.
Hence $B_{i,\nu}\cap B_{j,\nu}=\varnothing$.
Next, using \eqref{eq:centers}, we have 
\[
 |z_k(T_\nu)|=\sqrt{1-T_\nu}\,|z_k^0|
 =\nu^{1/2-\gmm}|z_k^0|.
\]
Since $C>\max_j|z_j^0|$ and $\gmm>0$, take $\nu$ smaller
if necessary so that
\[
 \nu^{\gmm/2}<C-\max_{1\le j\le3}|z_j^0|.
\]
For $x\in B_{k,\nu}$, \eqref{eq:observation-radius} then gives
\[
 \begin{aligned}
 |x|&\le |x-z_k(T_\nu)|+|z_k(T_\nu)|\\
 &<\nu^{(1-\gmm)/2}+\nu^{1/2-\gmm}|z_k^0|
 =\nu^{1/2-\gmm}\bigl(\nu^{\gmm/2}+|z_k^0|\bigr)\\
 &<C\nu^{1/2-\gmm}=R_\nu.
 \end{aligned}
\]
Thus $B_{k,\nu}\subset B(0,R_\nu)$ for every $k$.
We now estimate the Gaussian norm in \eqref{eq:gaussian-lower-bound}.
The decomposition, the disjointness of the balls, and
$B_{k,\nu}\subset B(0,R_\nu)$ give
\begin{equation}\label{eq:localization-reduction}
 \begin{aligned}
 \left\|\sum_k f_{k,\nu}\right\|_{L^1(B(0,R_\nu))}
 &=\left\|\sum_k f_{k,\nu}\mathbf1_{B_{k,\nu}}
 +\sum_k f_{k,\nu}\mathbf1_{\R^2\setminus B_{k,\nu}}
 \right\|_{L^1(B(0,R_\nu))}\\
 &\ge\left\|\sum_k f_{k,\nu}\mathbf1_{B_{k,\nu}}\right\|_{L^1(B(0,R_\nu))}
 -\left\|\sum_k f_{k,\nu}\mathbf1_{\R^2\setminus B_{k,\nu}}\right\|_{L^1(B(0,R_\nu))}\\
 &=\sum_k\int_{B_{k,\nu}}|f_{k,\nu}(x)|\dd x
 -\left\|\sum_k f_{k,\nu}\mathbf1_{\R^2\setminus B_{k,\nu}}\right\|_{L^1(B(0,R_\nu))}\\
 &\ge\sum_k\int_{B_{k,\nu}}|f_{k,\nu}(x)|\dd x
 -\sum_k\int_{\R^2\setminus B_{k,\nu}}|f_{k,\nu}(x)|\dd x\\
 &=\sum_k\|f_{k,\nu}\|_{L^1(\R^2)}
 -2\sum_k\|f_{k,\nu}\mathbf1_{\R^2\setminus B_{k,\nu}}\|_{L^1(\R^2)}.
 \end{aligned}
\end{equation}
The first sum in the last line is bounded below by
\eqref{eq:pair-limit}. We now estimate the second sum.
We first show that
\begin{equation}\label{eq:tail-domains}
 \begin{aligned}
 \R^2\setminus B_{k,\nu}
 &=\{x\in\R^2:|x-z_k(T_\nu)|\ge\nu^{(1-\gmm)/2}\}\\
 &\subset\{x\in\R^2:|x-\widehat z_k(T_\nu)|
 \ge\tfrac12\nu^{(1-\gmm)/2}\}.
 \end{aligned}
\end{equation}
Indeed, for $x\in\R^2\setminus B_{k,\nu}$,
\[
 \begin{aligned}
 |x-\widehat z_k(T_\nu)|
 &\ge |x-z_k(T_\nu)|-|\widehat z_k(T_\nu)-z_k(T_\nu)|\\
 &\overset{\eqref{eq:center-shift}}{\ge}
 \nu^{(1-\gmm)/2}-C_{\epsz}\sqrt\nu
 =(1-C_{\epsz}\nu^{\gmm/2})\nu^{(1-\gmm)/2}
 \ge\frac12\nu^{(1-\gmm)/2}.
 \end{aligned}
\]
By \eqref{eq:tail-domains}, the substitutions
$y=x-z_k(T_\nu)$ and $y=x-\widehat z_k(T_\nu)$, respectively, give
\[
 \begin{aligned}
 &\int_{\R^2\setminus B_{k,\nu}}
 H_{T_\nu}^\nu(x-z_k(T_\nu))\dd x
 =\int_{|y|\ge\nu^{(1-\gmm)/2}}H_{T_\nu}^\nu(y)\dd y,\\
 &\int_{\R^2\setminus B_{k,\nu}}
 H_{T_\nu-h_\nu}^{\nu/\lambda_\nu}(x-\widehat z_k(T_\nu))\dd x
 \le\int_{|y|\ge\frac12\nu^{(1-\gmm)/2}}
 H_{T_\nu-h_\nu}^{\nu/\lambda_\nu}(y)\dd y.
 \end{aligned}
\]
We now use the elementary Gaussian tail identity
\begin{equation}\label{eq:gaussian-tail}
 \int_{|x|\ge R}H_t^\mu(x)\dd x
 =e^{-R^2/(4\mu t)},\qquad \mu,t,R>0.
\end{equation}
It follows that
\begin{equation}\label{eq:pair-tails}
 \begin{aligned}
 &\sum_k\|f_{k,\nu}\mathbf1_{\R^2\setminus B_{k,\nu}}\|_{L^1(\R^2)}\\
 &\quad\overset{\eqref{eq:gaussian-pair},\,\eqref{eq:tail-domains}}{\le}
 \sum_k|\Gamma_k|\Biggl[
 \int_{|x|\ge\nu^{(1-\gmm)/2}}H_{T_\nu}^\nu(x)\dd x
 +\lambda_\nu\int_{|x|\ge\frac12\nu^{(1-\gmm)/2}}
 H_{T_\nu-h_\nu}^{\nu/\lambda_\nu}(x)\dd x\Biggr]\\
 &\quad\overset{\eqref{eq:gaussian-tail}}{=}\Gabs\left[
 e^{-\nu^{-\gmm}/(4T_\nu)}
 +\lambda_\nu e^{-\lambda_\nu\nu^{-\gmm}/(16(T_\nu-h_\nu))}
 \right]\\
 &\quad\le
 2\Gabs e^{-\nu^{-\gmm}/32}.
 \end{aligned}
\end{equation}
Combining \eqref{eq:gaussian-lower-bound},
\eqref{eq:localization-reduction}, \eqref{eq:pair-tails}, and
\eqref{eq:pair-limit}, we obtain
\[
 \|\omhat^\nu(T_\nu)
       -\omcol^\nu(T_\nu)\|_{L^1(B(0,R_\nu))}
 \ge2c_0-o(1)\ge c_0
\]
for sufficiently small $\nu$, uniformly for $t_0\in I$, where
\begin{equation}\label{eq:separation-constant}
 c_0=\sum_{k=1}^3|\Gamma_k|
 \erf\!\left(\frac{\epsz(1-\max I)|b_k|}{4}\right)>0.
\end{equation}
Here $c_0>0$ because at least one $b_k$ is nonzero.
This proves \eqref{eq:local-separation}. 
\end{proof}

\begin{proof}[Proof of Theorem~\ref{thm:main}]
Let $I$ be as in the statement and fix an arbitrary $t_0\in I$.
Use the parameters in \eqref{eq:parameters}, with $\epsz$
chosen as in \eqref{eq:choice-epsilon}, and let $\omhat^\nu$
be given by \eqref{eq:intermediate}.
We first construct the corrected solution and verify
\eqref{eq:circulation}--\eqref{eq:mass-bound}.
For $y\in\R^2$, let $\omtil^\nu(\cdot,\cdot;y)$ be the global smooth
solution of \eqref{eq:NS} with data
\begin{equation}\label{eq:corrected-data}
 \omtil^\nu(t_0,x;y)
 =\omhat^\nu(t_0,x)
 +(1-\lambda_\nu)\Gtot G(x-y).
\end{equation}
Conservation of circulation, $\int_{\R^2}G=1$, and
\eqref{eq:intermediate-mass} give
\[
 \int_{\R^2}\omtil^\nu(t,x;y)\dd x
 =\int_{\R^2}\omtil^\nu(t_0,x;y)\dd x
 =\lambda_\nu\Gtot
 +(1-\lambda_\nu)\Gtot
 =\Gtot,\qquad t\ge t_0.
\]
The $L^1$ maximum principle and \eqref{eq:intermediate-mass} yield
\[
 \|\omtil^\nu(t;y)\|_{L^1(\R^2)}
 \le\|\omtil^\nu(t_0;y)\|_{L^1(\R^2)}
 \le\lambda_\nu\Gabs
 +(1-\lambda_\nu)|\Gtot|
 \le\Gabs.
\]
Together with \eqref{eq:background-mass}, these prove
\eqref{eq:circulation} and \eqref{eq:mass-bound}.

We next prove \eqref{eq:initial-bound}. By \eqref{eq:corrected-data},
\eqref{eq:comparison-initial-explicit}, and \eqref{eq:early},
for $\nu$ sufficiently small that $C_I\nu^\gmm\le1/2$,
\[
 \begin{aligned}
 \|\omtil^\nu(t_0;y)-\omcol^\nu(t_0)\|_{L^1(\R^2)}
 &\le\|\omhat^\nu(t_0)-\omcol^\nu(t_0)\|_{L^1(\R^2)}
 +\epsz|\Gtot|\nu^{1-\gmm}\\
 &\le C_I\nu+
 \epsz\bigl(3\Gabs+|\Gtot|\bigr)
 \nu^{1-\gmm}\\
 &\overset{\eqref{eq:choice-epsilon}}{\le}
 \frac12\nu^{1-\gmm}+\frac12\nu^{1-\gmm}
 =\nu^{1-\gmm},
 \end{aligned}
\]
It remains to prove \eqref{eq:final-bound}. Adding and subtracting
$\omhat^\nu(T_\nu)$, we obtain
\begin{equation}\label{eq:separation-reduction}
 \begin{aligned}
 &\|\omtil^\nu(T_\nu;y)-\omcol^\nu(T_\nu)\|_{L^1(B(0,R_\nu))}\\
 &\quad=\left\|
 \bigl(\omhat^\nu(T_\nu)-\omcol^\nu(T_\nu)\bigr)
 +\bigl(\omtil^\nu(T_\nu;y)-\omhat^\nu(T_\nu)\bigr)
 \right\|_{L^1(B(0,R_\nu))}\\
 &\quad\ge
 \|\omhat^\nu(T_\nu)-\omcol^\nu(T_\nu)\|_{L^1(B(0,R_\nu))}
 -\|\omtil^\nu(T_\nu;y)-\omhat^\nu(T_\nu)\|_{L^1(B(0,R_\nu))}\\
 &\quad\overset{\eqref{eq:local-separation}}{\ge}
 c_0-\|\omtil^\nu(T_\nu;y)-\omhat^\nu(T_\nu)\|_{L^1(B(0,R_\nu))}.
 \end{aligned}
\end{equation}
We will choose $y=y_\nu$ so that the last norm tends to zero.
For each $y\in\R^2$, let
$g^\nu(\cdot,\cdot;y)$ solve the heat equation with the Gaussian
correction in \eqref{eq:corrected-data} as initial data:
\[
 \begin{cases}
 \partial_t g^\nu(t,x;y)=\nu\Delta_x g^\nu(t,x;y),
 &t>t_0,\quad x\in\R^2,\\
 g^\nu(t_0,x;y)=(1-\lambda_\nu)\Gtot G(x-y),
 &x\in\R^2.
 \end{cases}
\]
The heat evolution preserves the sign and integral of the initial data,
so
\begin{equation}\label{eq:correction-mass}
 \|g^\nu(t;y)\|_{L^1(\R^2)}
 =(1-\lambda_\nu)|\Gtot|
 =\epsz|\Gtot|\nu^{1-\gmm}.
\end{equation}
Since $g^\nu(t,\cdot;y)$ is radial about $y$, its velocity is tangential
and $v^{g^\nu}\cdot\nabla g^\nu=0$. Hence $g^\nu$ also solves \eqref{eq:NS}.
In differential identities below, we suppress the arguments $(t,x;y)$;
all spatial derivatives and convolutions are taken with $y$ fixed.
To compare the corrected solution with the sum of
$\omhat^\nu$ and this independently evolving Gaussian, define
\[
 \eta^\nu(t,x;y)=\omtil^\nu(t,x;y)
 -\omhat^\nu(t,x)-g^\nu(t,x;y),
 \qquad \eta^\nu(t_0;y)=0.
\]
We estimate this interaction error directly in $L^1(\R^2)$.
Subtracting the equations for $\omhat^\nu$ and $g^\nu$
from that for $\omtil^\nu$ gives
\begin{equation}\label{eq:error-equation}
 \partial_t\eta^\nu+v^{\omtil^\nu}\cdot\nabla\eta^\nu
 -\nu\Delta\eta^\nu
 =-v^{\eta^\nu}\cdot\nabla(\omhat^\nu+g^\nu)-R^\nu,
\end{equation}
where
\begin{equation}\label{eq:residual}
 R^\nu=v^{\omhat^\nu}\cdot\nabla g^\nu
 +v^{g^\nu}\cdot\nabla\omhat^\nu.
\end{equation}
Kato's inequality and $\nabla\cdot v^{\omtil^\nu}=0$ give,
for almost every $t$,
\begin{equation}\label{eq:error-energy}
 \frac{\dd}{\dd t}\|\eta^\nu\|_{L^1(\R^2)}
 \le\int_{\R^2}|v^{\eta^\nu}(x)|
 |\nabla(\omhat^\nu+g^\nu)(x)|\dd x
 +\|R^\nu\|_{L^1(\R^2)}.
\end{equation}
By Fubini's theorem and the Biot--Savart potential estimate, we obtain
\begin{equation}\label{eq:biot-bound}
 \begin{aligned}
 &\int_{\R^2}|v^{\eta^\nu}(x)|
 |\nabla(\omhat^\nu+g^\nu)(x)|\dd x\\
 &\quad\le\|\eta^\nu\|_{L^1(\R^2)}
 \bigl\||K|*|\nabla(\omhat^\nu+g^\nu)|\bigr\|_{L^\infty(\R^2)}\\
 &\quad\lesssim
 \|\nabla(\omhat^\nu+g^\nu)\|_{L^1(\R^2)}^{1/2}
 \|\nabla(\omhat^\nu+g^\nu)\|_{L^\infty(\R^2)}^{1/2}
 \|\eta^\nu\|_{L^1(\R^2)}.
 \end{aligned}
\end{equation}
For fixed $\nu,t_0$, translation invariance and contraction of the heat
semigroup in $W^{1,1}(\R^2)\cap W^{1,\infty}(\R^2)$, together with
the positive-time bounds for $\omhat^\nu$, which is independent
of $y$, give
\begin{equation}\label{eq:correction-gradient}
 \sup_{t_0\le t\le T_\nu}
 \left(
 \|\nabla\omhat^\nu(t)+\nabla g^\nu(t;y)\|_{L^1(\R^2)}
 +\|\nabla\omhat^\nu(t)+\nabla g^\nu(t;y)\|_{L^\infty(\R^2)}
 \right)\lesssim_{\nu,t_0}1,
\end{equation}
independently of $y$.
Since $\eta^\nu(t_0;y)=0$, \eqref{eq:error-energy}--\eqref{eq:correction-gradient}
and Gronwall's inequality yield
\begin{equation}\label{eq:comparison}
 \sup_{t_0\le t\le T_\nu}\|\eta^\nu(t;y)\|_{L^1(\R^2)}
 \le C_{\nu,t_0}\int_{t_0}^{T_\nu}\|R^\nu(t;y)\|_{L^1(\R^2)}\dd t.
\end{equation}
Here $C_{\nu,t_0}$ is finite for fixed $\nu,t_0$ and independent of $y$.
We now show that the integral in \eqref{eq:comparison} tends to zero
as $|y|\to\infty$, with $\nu,t_0$ fixed. Since
$g^\nu(t,x;y)=g^\nu(t,x-y;0)$, we have
$\nabla_xg^\nu(t,x;y)=\nabla_xg^\nu(t,x-y;0)$. The change of variables
$z=\xi+y$ also gives
\[
 \begin{aligned}
 v^{g^\nu}(t,x;y)
 &=\int_{\R^2}K(x-z)g^\nu(t,z-y;0)\dd z\\
 &=\int_{\R^2}K(x-y-\xi)g^\nu(t,\xi;0)\dd\xi
 =v^{g^\nu}(t,x-y;0).
 \end{aligned}
\]
Thus \eqref{eq:residual} and the change of variables $x=z+y$
in the first integral yield
\[
 \begin{aligned}
 \|R^\nu(t;y)\|_{L^1(\R^2)}
 &\le\int_{\R^2}|v^{\omhat^\nu}(t,x)|
 |\nabla g^\nu(t,x-y;0)|\dd x
 +\int_{\R^2}|v^{g^\nu}(t,x-y;0)|
 |\nabla\omhat^\nu(t,x)|\dd x\\
 &=\int_{\R^2}|v^{\omhat^\nu}(t,z+y)|
 |\nabla g^\nu(t,z;0)|\dd z
 +\int_{\R^2}|v^{g^\nu}(t,x-y;0)|
 |\nabla\omhat^\nu(t,x)|\dd x.
 \end{aligned}
\]
Fix $t\in[t_0,T_\nu]$. Both velocities are bounded and vanish at
spatial infinity, so for each fixed $z,x\in\R^2$,
\[
 v^{\omhat^\nu}(t,z+y)\longrightarrow0,
 \qquad v^{g^\nu}(t,x-y;0)\longrightarrow0
 \qquad (|y|\to\infty).
\]
Moreover,
\[
 \begin{aligned}
 |v^{\omhat^\nu}(t,z+y)|\,|\nabla g^\nu(t,z;0)|
 &\le\|v^{\omhat^\nu}(t)\|_{L^\infty(\R^2)}
 |\nabla g^\nu(t,z;0)|,\\
 |v^{g^\nu}(t,x-y;0)|\,|\nabla\omhat^\nu(t,x)|
 &\le\|v^{g^\nu}(t;0)\|_{L^\infty(\R^2)}
 |\nabla\omhat^\nu(t,x)|.
 \end{aligned}
\]
The right-hand sides are independent of $y$ and integrable in $z$
and $x$, respectively. Dominated convergence in each integral gives
\[
 \|R^\nu(t;y)\|_{L^1(\R^2)}\longrightarrow0
 \qquad (|y|\to\infty).
\]
To pass to the time integral, H\"older's inequality gives
\[
 \begin{aligned}
 \|R^\nu(t;y)\|_{L^1(\R^2)}
 &\le\|v^{\omhat^\nu}(t)\|_{L^\infty(\R^2)}
 \|\nabla g^\nu(t;0)\|_{L^1(\R^2)}\\
 &\quad+\|v^{g^\nu}(t;0)\|_{L^\infty(\R^2)}
 \|\nabla\omhat^\nu(t)\|_{L^1(\R^2)}
 \lesssim_{\nu,t_0}1.
 \end{aligned}
\]
Indeed, the Biot--Savart estimate
$\|K*f\|_{L^\infty(\R^2)}\lesssim
\|f\|_{L^1(\R^2)}^{1/2}\|f\|_{L^\infty(\R^2)}^{1/2}$,
the bounds for the reference solution at positive times, and
contraction of the heat semigroup bound all four factors uniformly
in $t\in[t_0,T_\nu]$. This constant is independent of $y$ and integrable
on the finite time interval. Dominated convergence in time therefore yields
\[
 \int_{t_0}^{T_\nu}\|R^\nu(t;y)\|_{L^1(\R^2)}\dd t
 \longrightarrow0\qquad (|y|\to\infty).
\]
Substituting this limit into \eqref{eq:comparison}, we obtain
\begin{equation}\label{eq:decoupling}
 \begin{aligned}
 &\sup_{t_0\le t\le T_\nu}
 \|\omtil^\nu(t;y)-\omhat^\nu(t)-g^\nu(t;y)\|_{L^1(\R^2)}\\
 &\quad\le C_{\nu,t_0}
 \int_{t_0}^{T_\nu}\|R^\nu(t;y)\|_{L^1(\R^2)}\dd t
 \longrightarrow0\qquad (|y|\to\infty).
 \end{aligned}
\end{equation}
More precisely, for each fixed $\nu$ and $t_0$, choose
$y=y_\nu(t_0)$ sufficiently far away that
\begin{equation}\label{eq:interaction-choice}
 \int_{t_0}^{T_\nu}\|R^\nu(t;y_\nu)\|_{L^1(\R^2)}\dd t
 \le\frac{\nu}{1+C_{\nu,t_0}},
\end{equation}
where $C_{\nu,t_0}$ is the constant in \eqref{eq:comparison}.
Using \eqref{eq:comparison} and \eqref{eq:interaction-choice}, we obtain
\begin{equation}\label{eq:error-choice}
 \sup_{t_0\le t\le T_\nu}\|\eta^\nu(t;y_\nu)\|_{L^1(\R^2)}
 \le C_{\nu,t_0}\frac{\nu}{1+C_{\nu,t_0}}\le\nu.
\end{equation}
Moving the correction farther away preserves its circulation and
$L^1(\R^2)$ norm, so the initial bound remains valid.
Since $B(0,R_\nu)\subset\R^2$, the triangle inequality gives
\begin{equation}\label{eq:local-correction}
 \begin{aligned}
 &\|\omtil^\nu(T_\nu;y_\nu)-\omhat^\nu(T_\nu)\|_{L^1(B(0,R_\nu))}\\
 &\quad\le\|\eta^\nu(T_\nu;y_\nu)\|_{L^1(\R^2)}
 +\|g^\nu(T_\nu;y_\nu)\|_{L^1(\R^2)}\\
 &\quad\overset{\eqref{eq:error-choice},\,\eqref{eq:correction-mass}}{\le}
 \nu+\epsz|\Gtot|\nu^{1-\gmm}
 \longrightarrow0.
 \end{aligned}
\end{equation}
Taking $\omtil^\nu(t,x)=\omtil^\nu(t,x;y_\nu)$,
\eqref{eq:separation-reduction} and \eqref{eq:local-correction} yield
\[
 \|\omtil^\nu(T_\nu)-\omcol^\nu(T_\nu)\|_{L^1(B(0,R_\nu))}
 \ge c_0-\nu-\epsz|\Gtot|\nu^{1-\gmm}
 \ge\frac{c_0}{2}
\]
for sufficiently small $\nu$, uniformly for $t_0\in I$.
Thus \eqref{eq:final-bound} holds with $c=c_0/2$, where $c_0$
is given by \eqref{eq:separation-constant}.
The constant $c$ is independent of $\gmm$.
All smallness conditions on $\nu$ are uniform for $t_0\in I$.
This completes the proof.
\end{proof}

\subsection*{Acknowledgments}
The author is supported by the NSFC Excellent Young Scientists
Fund (Overseas) and the Research Start-up Fund of Shanghai Jiao Tong
University. 

\subsubsection*{Data availability}
No datasets were generated or analysed during the current study.

\subsubsection*{Competing interests}
The author declares no competing interests.

\end{document}